\documentclass[conference]{IEEEtran}
\usepackage{amsmath,graphicx}
\usepackage{amsfonts,amssymb}
\usepackage{multirow}
\usepackage{algorithm}
\usepackage{algpseudocode}
\usepackage{cite}
\usepackage{balance}
\DeclareGraphicsExtensions{.jpg,pdf,.png}
\ifCLASSINFOpdf
\else
\fi
\usepackage{fancyhdr}
\newcommand*{\myfont}{\fontfamily{phv}\selectfont}
\DeclareTextFontCommand{\textmyfont}{\myfont}

\begin{document}
\title{Adaptive Singular Value Thresholding}
\author{Nematollah Zarmehi and Farokh Marvasti\\Advanced Communication Research Institute\\Department of Electrical Engineering\\Sharif University of Technology, Tehran, Iran\\Email: http://zarmehi.ir/contact.html, marvasti@sharif.edu}
\maketitle
\thispagestyle{fancy}
\begin{abstract}
In this paper, we propose an Adaptive Singular Value Thresholding (ASVT) for low rank recovery under affine constraints. Unlike previous iterative methods that the threshold level is independent of the iteration number, in our proposed method, the threshold in adaptively decreases during iterations. The simulation results reveal that we get better performance with this thresholding strategy.
\end{abstract}
\section{Introduction}\label{sec:intro}
Finding a minimum rank matrix subject to affine constraints is the heart of many engineering applications such as minimal realization theory \cite{ref:minreal}, minimum order controller design \cite{ref:mincont1,ref:mincont2,ref:mincont3}, collaborative filtering \cite{ref:collfilt}, and machine learning \cite{ref:machl1,ref:machl2,ref:machl3}. This problem is usually known as Affine Rank Minimization (ARM). Mathematically speaking, we want to solve the following minimization problem:
\begin{equation}\label{eq:arm}
\begin{split}
\min_\mathbf{X} & ~ rank(\mathbf{X}), \\
s.t.: & ~ \mathbf{\mathcal{A}}(\mathbf{X})=\mathbf{b},
\end{split}
\end{equation}
where $\mathbf{X}\in \mathbb{R}^{n_1\times n_2}$ is the unknown matrix, $\mathbf{\mathcal{A}}:\mathbb{R}^{n_1\times n_2} \to \mathbb{R}^m$ denotes the linear mapping, and $\mathbf{b}$ is the measurement vector. The affine constraint $\mathbf{\mathcal{A}}(\mathbf{X})=\mathbf{b}$ can be simplified as
\begin{equation*}
\mathbf{A}vec(\mathbf{X})=\mathbf{v}
\end{equation*}
where $\mathbf{A} \in \mathbb{R}^{m \times n_1n_2}$ is the matrix representation of $\mathbf{\mathcal{A}}$ and $vec(.)$ is the standard vectorization operator. 

Since rank is the number of non-zero Singular Values (SVs), the minimization problem (\ref{eq:arm}) is the same as
\begin{equation}\label{eq:arm2}
\begin{split}
\min_\mathbf{X} & ~ \|\mathbf{\sigma}(\mathbf{X})\|_0, \\
s.t.: & ~ \mathbf{\mathcal{A}}(\mathbf{X})=\mathbf{b},
\end{split}
\end{equation}
where $\mathbf{(\sigma)}(\mathbf{X})$ is the vector of SVs of $\mathbf{X}$.

An important special case of ARM problem introduced in (\ref{eq:arm}) is the matrix completion problem:
\begin{equation}\label{eq:mc}
\begin{split}
\min_\mathbf{X} & ~ rank(\mathbf{X}), \\
s.t.: & ~ [\mathbf{X}]_{i,j}=[\mathbf{X_*}]_{i,j}, \quad \forall(i,j)\in \Omega,
\end{split}
\end{equation}
where $\mathbf{X}$ is the unknown matrix, $[\mathbf{X}]_{i,j}$ is the element of matrix $\mathbf{X}$ at the intersection of $i$-th row and $j$-th column, and $\mathbf{X_*}$ is the low rank matrix that we would like to recover but we only know a small subset of all entries of $\mathbf{X_*}$, that is, $\Omega \subseteq \{1,2,\cdots,n_1\}\times \{1,2,\cdots,n_2\}$. 

There are about two main approaches to solve the ARM problem or better saying approximate the non-convex minimization problem (\ref{eq:arm2}).

In the first approach, the rank function is approximated by the nuclear norm. The nuclear norm of a matrix is defined as $\|\mathbf{X}\|_*=\sum_{i=1}^{r}\sigma_i(\mathbf{X})$ where $\sigma_i$ is the $i$-th largest SV of $\mathbf{X}$. Therefore, we can write the nuclear norm as $\|\mathbf{X}\|_*=\|\sigma(\mathbf{X})\|_1$ and finally, the nuclear norm minimization will be equivalent to an $\ell_1$-norm minimization as follows \cite{ref:fazelboyd,ref:mfazel}:
\begin{equation}\label{eq:l1}
\begin{split}
\min_\mathbf{X} & ~ \|\mathbf{\sigma}(\mathbf{X})\|_1, \\
s.t.: & ~ \mathbf{\mathcal{A}}(\mathbf{X})=\mathbf{b}.
\end{split}
\end{equation}
Above nuclear norm minimization can be efficiently solved using some convex optimization techniques such as a semi-definite program \cite{ref:fazel2}. For the purpose of matrix completion, a Singular Value Thresholding (SVT) algorithm is proposed to solve the nuclear norm minimization \cite{ref:svt}. It is an iterative algorithm that for a fixed value of $\tau>0$ and positive step-sizes $\{\delta_k\}$, starts with $\mathbf{Y}_0=\mathbf{0}$ and repeats
\begin{equation}
\left\{ \begin{array}{l}
\mathbf{X}_k = \mathcal{D}(\mathbf{Y}_{k-1},\tau)\\
~\\
\mathbf{Y}_k  = \mathbf{Y}_{k-1}+\delta_k\mathbf{\mathcal{A}}(\mathbf{X_*}-\mathbf{X}_k)
\end{array} \right.
\end{equation}
until a stopping criterion is satisfied. $\mathcal{D}(.,.)$ is called the singular value shrinkage operator. For a matrix $\mathbf{A} \in \mathbb{R}^{n_1\times n_2}$ with Singular Value Decomposition (SVD) of $\mathbf{A}=\mathbf{U}\mathbf{\Sigma}\mathbf{V}^T$, we have \begin{equation}\label{eq:shrnkg}
\mathcal{D}(\mathbf{A},\tau) = \mathbf{U}diag(\{\max(\sigma_i-\tau,0)\}_{i=1}^{\min(n_1,n_2)})\mathbf{V}^T.
\end{equation}

In the second approach, the $\ell_0$-norm used in the rank function is approximated by families of smoothing $\ell_0$-norm functions \cite{ref:malek}. For example, a class of Gaussian smoothing functions are defined as:
\begin{equation}
f_\delta (x) = \exp\left({-\frac{x^2}{2\delta^2}}\right).
\end{equation}
As $\delta \to 0$, $f_\delta(.)$ gives better approximation of the Kronecker delta function. In this approach, the rank of $\mathbf{X}\in \mathbb{R}^{n_1\times n_2}$ is approximated as follows:
\begin{equation}
rank(\mathbf{X}) \approx min(n_1,n_2)-\sum_{i=1}^{min(n_1,n_2)}f_\delta\big(\sigma_i(\mathbf{X})\big)
\end{equation}
and the final optimization problem is solved using the Gradient Projection (GP) algorithm \cite{ref:gp}. 

In this paper, we consider the ARM problem (\ref{eq:arm}) and propose an adaptive singular value thresholding method. Here, the threshold level, unlike the threshold level in the SVT algorithm is not fixed. 

\section{The Proposed Algorithm}\label{sec:proposed}
We aim to recover a low rank matrix under some linear constraints. As stated in Section \ref{sec:intro}, this problem is known as ARM. It should be noted that we don't know the rank of the primary matrix. Since the low rank matrix has sparsity in SV domain, we propose to threshold the SVs. But the main difference between our algorithm and previous ones is that the threshold level is adaptively changed during iterations. The idea of adaptive thresholding is used in \cite{ref:imat} for recovery of the original signal from its random samples. As suggested in \cite{ref:imat}, we choose threshold level of $k$-th iteration as follows:
\begin{equation}
\tau_k = b\exp(-ak),
\end{equation}
where $a$ and $b$ are two constants. The proposed algorithm is as follows:
\begin{equation}\label{eq:proposed}
\left\{ \begin{array}{l}
\mathbf{X}_k = \mathcal{T}(\mathbf{Y}_{k-1},\tau_k)\\
~\\
\mathbf{Y}_k  = \mathbf{Y}_{k-1}+\delta_k\mathbf{\mathcal{A}}^*(\mathbf{b}-\mathbf{\mathcal{A}}\mathbf{X}_k)
\end{array} \right.
\end{equation}
where $\mathbf{\mathcal{A}}^*(.)$ is the adjoint of $\mathbf{\mathcal{A}}(.)$. For a matrix $\mathbf{X}\in \mathbb{R}^{n_1\times n_2}$ with SVD $\mathbf{X}=\mathbf{U}\mathbf{\Sigma}\mathbf{V}^T$, $\mathcal{T}(\mathbf{X},\tau_k)$ is defined as follows:
\begin{equation}
\mathcal{T}(\mathbf{X},\tau_k) \buildrel \Delta \over = \mathbf{U}\big[\mbox{if~}\sigma_{i,j}<\tau_k: 0, \mbox{~else:~}\sigma_{i,j}\big]_{i=1,j=1}^{n_1,n_2}\mathbf{V}^T.
\end{equation}

We call our algorithm Adaptive Singular Value Thresholding (ASVT) which is shown in Algorithm 1. 

{\linespread{1.4}
\begin{algorithm}\label{alg}
\caption{Adaptive Singular Value Thresholding: ASVT}
\begin{algorithmic}[1]
\State \textbf{input:}
\State \indent {Measurement vector:} $\mathbf{b} \in \mathbb{R}^{m}$
\State \indent {Linear mapping:} $\mathbf{\mathcal{A}}$
\State \indent {Maximum number of iterations:} $K$
\State \indent {Positive step-sizes:} $\{\delta_k\}_{k=1}^{K}$
\State \indent {Two constants for thresholding:} $a, b$
\State \indent {Termination tolerance:} $\epsilon$
\State \textbf{initialization:}
\State \indent $k\gets 1$
\State \indent $e\gets \infty$
\State \indent $\mathbf{X}_0 \gets \mathbf{0}$
\State \indent $\mathbf{Y} \gets \mathbf{0}$
\While {$e>\epsilon ~\&~ k<K$}
\State $[\mathbf{U},\mathbf{\Sigma},\mathbf{V}]=svd(\mathbf{Y})$
\State $\tau_k=b\exp(-ak)$
\State $\mathbf{\Sigma}(\mathbf{\Sigma}<\tau_k) = 0$
\State $\mathbf{X}_k = \mathbf{U}\mathbf{\Sigma}\mathbf{V}^T$
\State $\mathbf{Y}=\mathbf{Y}+\delta_k\mathbf{\mathcal{A}}^*(\mathbf{b}-\mathbf{\mathcal{A}}\mathbf{X}_k)$
\State $e=\|\mathbf{X}_k-\mathbf{X}_{k-1}\|_F$
\State $k\gets k+1$
\EndWhile
\State \indent $\mathbf{\hat{X}}\gets \mathbf{X}_{k-1}$
\State \textbf{return} $\mathbf{\hat{X}}$
\end{algorithmic}
\end{algorithm}}
	
\section{Simulation Results}
The numerical experiment results are presented in this section. All simulations are done by MATLAB R2015a on Intel(R) Core(TM) i7-5960X @ 3GHz with 32GB-RAM.
	
We consider matrix completion as a special case of linear constraints. The random matrices are generated according to the scheme proposed in \cite{ref:svt}. An $r$-rank matrix $\mathbf{X}$ of size $n_1\times n_2$ is generated as $\mathbf{X}=\mathbf{A}\mathbf{B}$ where $\mathbf{A}\in \mathbb{R}^{n_1\times r}$ and $\mathbf{B}\in \mathbb{R}^{r\times n_2}$ are sampled independently form a standard Gaussian distribution ($\mathcal{N}(0,1)$). The measurements or the subset of observed entries $\Omega$ is sampled uniformly at random of size $m$. As \cite{ref:svt}, we choose constant step-sizes as $\delta = 1$.

To evaluate the performance of the algorithms, we report the relative error of the reconstruction defined as follows:
\begin{equation}
RE(\mathbf{X},\mathbf{\hat{X}}) = \frac{\|\mathbf{\hat{X}}-\mathbf{X}\|_F}{\|\mathbf{X}\|_F}
\end{equation}
where $\mathbf{X}$ and $\mathbf{\hat{X}}$ are original and reconstructed matrices, respectively.
	
\subsection{Comparison with the SVT Algorithm}\label{subsec:comp}
We generate matrices of size $n_1 = n_2 = 100, 300, 500, 1000, \mbox{and},2000$ with different ranks and different fraction of observations $\left({\frac{m}{n_1n_2}}\right)$ and provide a comparison between our algorithm with the SVT algorithm proposed in \cite{ref:svt}. We have downloaded the SVT MATLAB codes from http://svt.stanford.edu/code.html.

\begin{table}[h!]
\renewcommand{\arraystretch}{1.3}
\centering
\caption{Comparison of the proposed algorithm and the SVT algorithm.}\label{tab:comp1}\vspace{0.2cm}
\begin{tabular}{|c|c|c|c|c||c|c|}
\cline{4-7}
\multicolumn{3}{c|}{} & \multicolumn{2}{|c||}{ASVT} & \multicolumn{2}{|c|}{SVT \cite{ref:svt}} \\ \hline 
\scriptsize{\textbf{Size}} & \scriptsize{\textbf{Rank}} & \scriptsize{$\frac{m}{n_1n_2}$} & \scriptsize{$\#$iters} & \scriptsize{RE} & \scriptsize{$\#$iters} & \scriptsize{RE} \\ \hline \hline
\multirow{3}{*}{$500\times 500$}& 10 & 0.15 & 75 & 9.60e-4 & 200 & 5.02e-3\\ 
&  50 & 0.4 & 64 & 9.62e-4 & 200 & 4.38e-2\\
& 100 & 0.5 & 138 & 9.64e-4 & 200 & 3.28e-1\\ \hline
\multirow{3}{*}{$1000\times 1000$}& 10 & 0.15 & 38 & 9.22e-4 & 62 & 1.49e-3\\ 
&  50 & 0.4 & 29 & 9.21e-4 & 70 & 1.60e-3\\
& 100 & 0.5 & 37 & 9.88e-4 & 103 & 1.89e-3\\ \hline
\multirow{3}{*}{$1500\times 1500$}& 10 & 0.15 & 23 & 8.17e-4 & 49 & 1.28e-3\\ 
&  50 & 0.4 & 22 & 8.69e-4 & 53 & 1.32e-3\\
& 100 & 0.5 & 26 & 8.79e-4 & 68 & 1.51e-3\\ \hline
\multirow{3}{*}{$2000\times 2000$}& 10 & 0.1 & 49 & 9.24e-4 & 54 & 1.31e-3\\ 
&  50 & 0.3 & 27 & 8.05e-4 & 56 & 1.35e-3\\
& 100 & 0.4 & 28 & 9.98e-4 & 68 & 1.53e-3\\ \hline
\multirow{3}{*}{$2500\times 2500$}& 10 & 0.1 & 38 & 9.60e-4 & 47 & 1.65e-3\\ 
&  50 & 0.25 & 33 & 8.36e-4 & 56 & 1.37e-3\\
& 100 & 0.3 & 37 & 8.96e-4 & 77 & 1.61e-3\\ \hline
\multirow{3}{*}{$3000\times 3000$}& 10 & 0.05 & 29 & 9.29e-4 & 65 & 1.53e-3\\ 
&  50 & 0.25 & 29 & 9.33e-4 & 50 & 1.33e-3\\
& 100 & 0.3 & 32 & 8.37e-4 & 66 & 1.51e-3\\ \hline
\end{tabular}
\end{table}
	
According to the results of Table \ref{tab:comp1}, our proposed method has better performance in terms of relative error even with a fewer number of iterations. For example, when the size of the matrix is $3000\times 3000$ and its rank is 100 and only 0.3 of its entries are observed, our algorithm recovers it after 29 iterations while the SVT algorithm does this after 65 iterations. Moreover, the relative error of our algorithm and the SVT algorithm is $9.29\times 10^{-4}$ and $1.53\times 10^{-3}$, respectively. 
	
\subsection{Effect of Parameters}\label{subsec:param}
\begin{figure}[!h]
\centering
\includegraphics[width=1.02\linewidth]{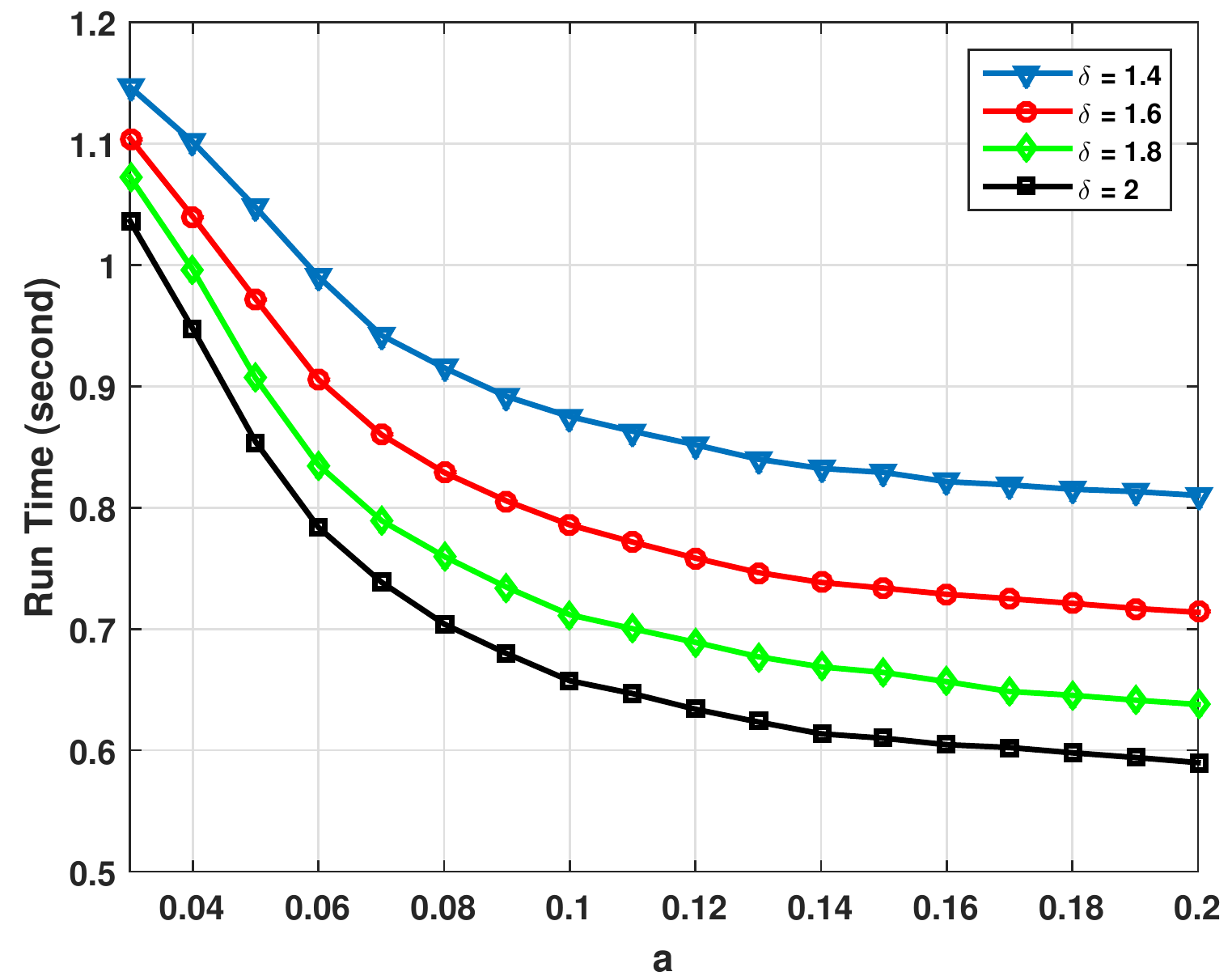}
\caption{Effect of damping factor $\alpha$.}\label{fig:effa}
\end{figure}

\begin{figure}[!h]
\centering
\includegraphics[width=1.02\linewidth]{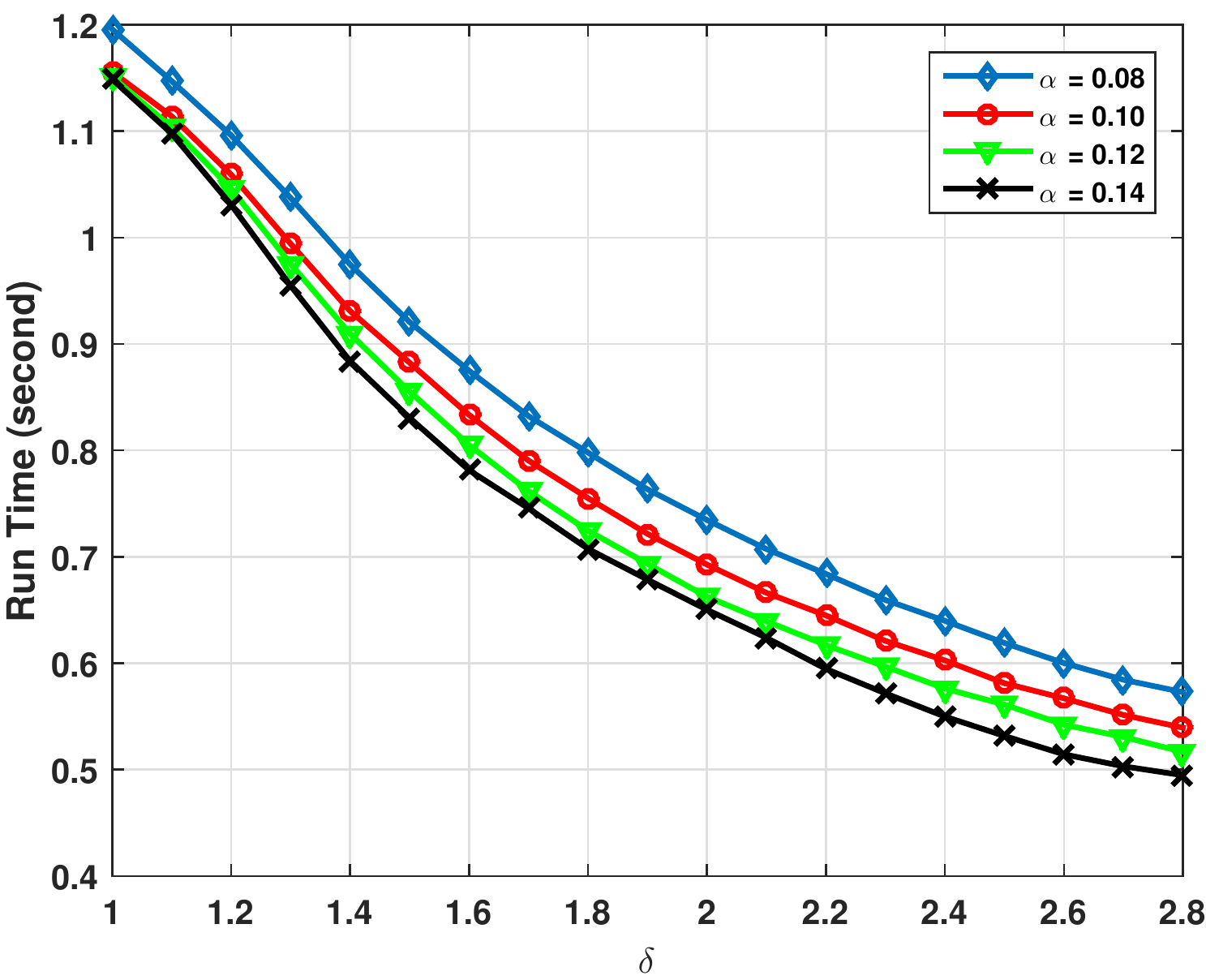}
\caption{Effect of step size $\delta$.}\label{fig:effdelta}
\end{figure}

In this Subsection, we investigate the effect of damping factor in thresholding operator ($\alpha$) and step size ($\delta$). For this simulation, we generate $300 \times 300$ random matrices of rank $10$ where only $1/3$ of their entries are observed. We test the proposed method with different values of $\alpha$ and $\delta$. The simulation results are shown in Figs. \ref{fig:effa} and \ref{fig:effdelta}. According to the results of Figs. \ref{fig:effa} and \ref{fig:effdelta}, the run time of the algorithm decrease as $\alpha$ and $\delta$ increases.
	
\subsection{Phase Transition Plot}\label{subsec:obs}
To better show the performance of the proposed method, we plot the phase transition between $\frac{m}{n_1n_2}$ and $\frac{d_r}{m}$ where $d_r=r(n_1+n_2-r)$ is the degree of freedom.

\begin{figure}[!h]
\centering
\includegraphics[width=1.02\linewidth]{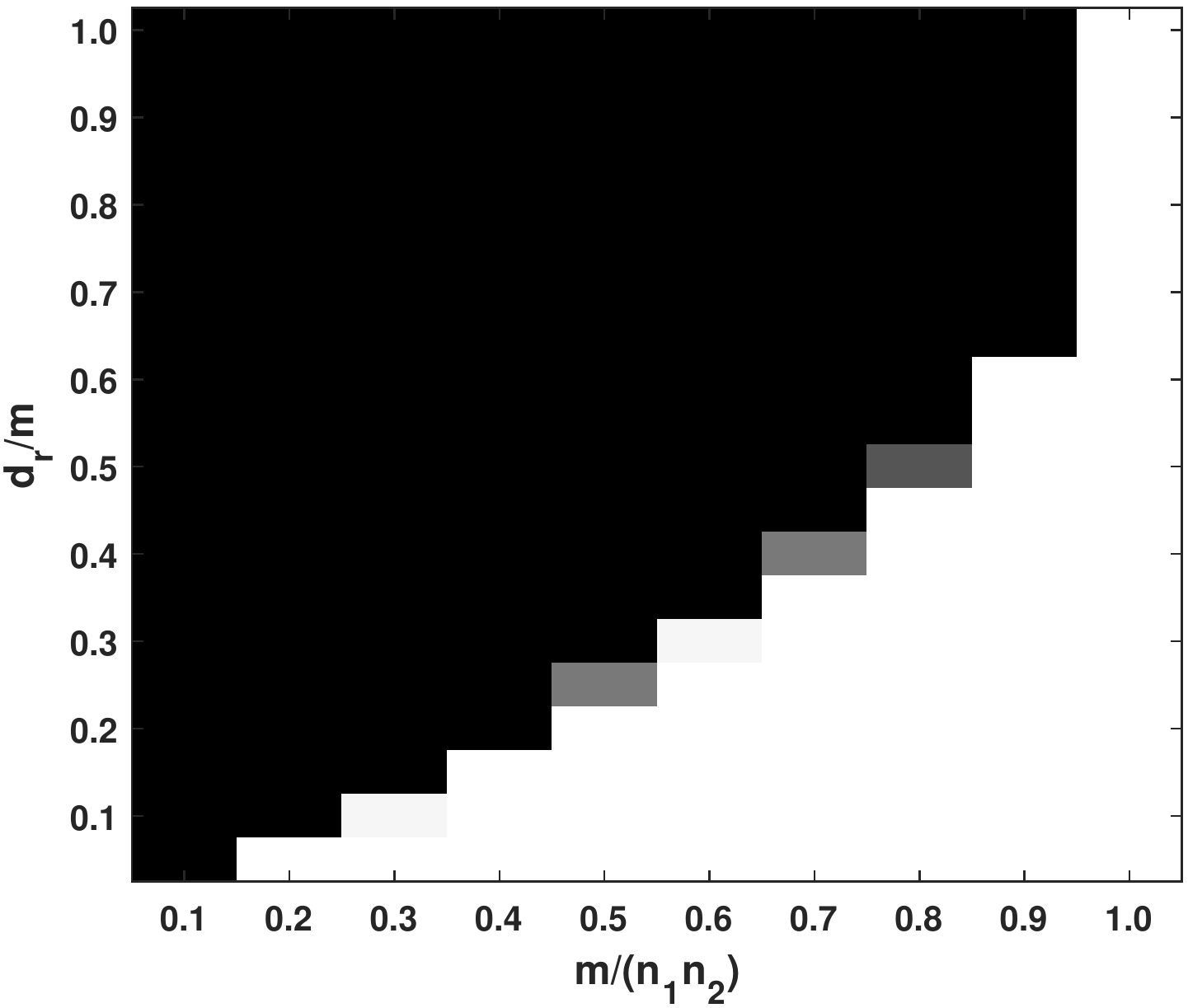}
\caption{Phase transition plot between $d_r/m$ and $m/(n_1n_2)$.}\label{fig:phtr}
\end{figure}

In this simulation, we generate random matrix $\mathbf{M}$ of size $80\times 80$. Each experiment is repeated 100 times and a matrix is declared to be successfully recovered if the relative error between the recovered matrix and the original one is less than $10^{-3}$. The phase transition is charted in Fig. \ref{fig:phtr}. The gray color of each cell indicates the probability of success recovery during simulations. For example, the white color indicates that the proposed algorithm can successfully recover the original matrix while the black indicates that it fails.
	
\section{Conclusion}
In this paper, we proposed an adaptive singular value thresholding for recovery of low rank matrices under affine constraints. The word adaptive means that the threshold level is adaptively changed during the iterations. Recently, the SVT algorithm is proposed for the matrix completion which is an important case of low rank matrix recovery under affine constraints. We also compared our algorithm with the SVT algorithm and the simulation results show that we can achieve improvement with adaptive thresholding. We show that the relative error of the proposed method is less than the relative error of the SVT algorithm. Moreover, our algorithm does this task in a fewer number of iterations.
	

\end{document}